\documentclass[12pt]{article} 
\usepackage{amsmath,amsthm, amssymb} 
\usepackage{amssymb,latexsym}

\newtheorem*{theorem}{Theorem}

\newtheorem{lemma}{Lemma}

\begin{document}

\baselineskip=17pt 

\title{\bf On the number of pairs of positive integers $x, y \le H$ such that
$x^2 + y^2 + 1$ is squarefree}

\author{\bf D. I. Tolev}

\date{}
\maketitle

\begin{abstract}
It is not difficult to find an asymptotic formula for the number of pairs of positive integers $x, y \le H$ such that
$x^2 + y^2 + 1$ is squarefree. In the present paper we improve the estimate for the error term in this formula
using the properties of certain exponential sums. A.Weils's estimate for the Kloosterman sum plays the major role in our analysis.

\smallskip


\end{abstract}

\section{Notations}

\indent

Let $H$ be a sufficiently large positive number. 
By the letters $k, m, n$ we denote integers and by $ d, r, l, h, q, x, y $ --- positive integers. 
The letters $t$ and $D$ denote real numbers and the letter $p$ is reserved for primes.
By $\varepsilon$ we denote an arbitrary small positive number, not necessarily the same in different occurrences. 
This convention allows us to write $q^{\varepsilon} \log q \ll q^{\varepsilon}$, for example.
If it is not explicitly stated the opposite, the constants in the Vinogradov and Landau symbols are absolute or depend on $\varepsilon$.

\bigskip

We denote by $\mu(n)$ the M\"obius function and by $\tau(n)$ the number of positive divisors of $n$. 
We write $(n_1, \dots, n_k)$ for the greatest common divisor of $n_1, \dots, n_k$.
As usual $[ t]$ and $\{ t \}$ denote the integer part, respectively, the fractional part of $t$.
We put $\rho(t) = \frac{1}{2} - \{ t \}$ and let 
$|| t ||$ be the distance from $t$ to the nearest integer.
Further $e(t) = \exp \left( 2 \pi i t \right) $ and $e_q(t) = e(t/q)$. 
For any $q$ and $k$ such that $(q, k)=1$ we denote by $\overline{(k)}_q$ 
the inverse of $k$ modulo $q$. If the value of the modulus is clear form the context then we write for simplicity $\overline{k}$. For any odd $q$ we denote by $\left( \frac{\cdot}{q}\right) $ the Jacobi symbol.
By $\square$ we mark an end of a proof or its absence.

\bigskip

In this paper we use the properties of the Gauss sum and the Kloosterman sum. 
They are defined by
\begin{equation} \label{100}
  G(q; n, m) = \sum_{1 \le x \le q} e_q \left( n x^2 + m x\right) , \qquad
  G(q; n) = G(q; n, 0) 
\end{equation}
and, respectively, by
\begin{equation} \label{200}
  K(q; n, m) = \sum_{\substack{ 1 \le x \le q \\ (x, q) = 1 }} e_q \left( n x + m \overline{x} \right) .
\end{equation}

\section{Introduction and statement of the result}

\indent

Suppose that $f(t_1, \dots, t_r)$ is a polynomial with integer coefficients
and denote by $S_f(H)$ the number of $r$-tuples of positive integers $n_1, \dots, n_r \le H$ such that 
$f(n_1, \dots, n_r)$ is squarefree. 
The problem is to evaluate $S_f(H)$ when $H$ is large.
A lot of articles are devoted to problems of this type.
We point out the papers of
Estermann~\cite{Est1}, Filaseta~\cite{Filaseta}, Greaves~\cite{Greaves}, Heath-Brown~\cite{Heath-Brown},
Hooley~\cite{Hooley1}, \cite{Hooley2}, (see also the book \cite{Hooley}, Chapter~4)  and Poonen~\cite{Poon}, 
but many other similar results can be found in literature.

\bigskip

In certain cases one can obtain an asymptotic formula for $S_f(H)$.
For example Estermann~\cite{Est1} proved that
\[
  \sum_{1 \le x \le H} \mu^2 \left( x^2 + 1 \right) = c_0 H + O \left( H^{ \frac{2}{3} + \varepsilon}\right) ,
\]
where $c_0>0$ is a constant. 
We should mention that using his ``square sieve'' Heath-Brown's~\cite{Heath-Brown}
considered $S_f(H)$ for the polynomial $f(t) = t(t+1)$ and found
an asymptotic formula with error term $ O \big( H^{\frac{7}{11} + \varepsilon }\big) $.

\bigskip

In the present paper we consider $S_f(H)$ for the polynomial
$f(t_1, t_2) = t_1^2 + t_2^2 + 1$. We write for simplicity
\[
  S(H) = \sum_{1 \le x, y \le H} \mu^2(x^2 + y^2 + 1) 
\]
and denote
\begin{equation} \label{300}
  \lambda(q; n, m) = \sum_{ x, \, y : \; \eqref{400}} 
   e_q \left( n x + m y \right) ,
\end{equation}
where the summation is taken over the integers $x, y$ satisfying the conditions
\begin{equation} \label{400}
1 \le x, y \le q , \qquad x^2 + y^2 + 1 \equiv 0 \pmod{q} .
\end{equation}
We denote also
\begin{equation} \label{405}
 \qquad
   \lambda(q) = \lambda(q; 0, 0) .
\end{equation}

\bigskip

Our result is the following
\begin{theorem}
For the sum
$S(H)$ 
we have 
\begin{equation} \label{410}
   S(H) = c H^2 + O \left( H^{\frac{4}{3} + \varepsilon } \right) ,
\end{equation}
where 
\begin{equation} \label{412}
  c = \prod_{p} \left( 1 - \frac{\lambda \left( p^2\right) }{p^4} \right) .
\end{equation}
\end{theorem}

\bigskip

We note that if we apply simplest elementary methods only then we will find an asymptotic formula 
for $S(H)$ with an error term $O \big( H^{\frac{3}{2} + \varepsilon} \big)$. 
Our better result is a consequence of the estimate for $\lambda(q; n, m)$ established in Lemma~\ref{est-lambda},
which can be obtained using the properties of the Gauss sum and A.Weil's estimate for the Kloosterman sum.

\bigskip

We may apply the same method for studying $S_f(H)$ 
with any quadratic polynomial $f$ in two variables. 
(Of course trivial cases like $f(t_1, t_2) = (t_1+t_2)^2$ have to be excluded.)
Then the exponential sum
\begin{equation} \label{414}
  \lambda_f (q; n, m) = \sum_{\substack{ 1 \le x, y \le q \\ f(x, y) \equiv 0 \pmod{q}} }
    e_q (nx + my) 
\end{equation}
naturally appears. 
It is closely connected with the Kloosterman sum (and in the case $f(t_1, t_2) = t_1 t_2 - 1 $ the sum \eqref{414} coincides with the Kloosterman sum),
therefore we may use again A.Weil's estimate in order to estimate \eqref{414}.
However we shall not consider this more general problem here.

\section{Lemmas}

\indent

Our first lemma includes the basic properties of the Gauss sum. The proofs are available in  Section~6 of \cite{Est2}
and Chapter~7 of \cite{Hua}.
\begin{lemma} \label{Gauss}
For the Gauss sum we have

(i)
If $(q, n)=d$ then
\[
  G(q; n, m) = \begin{cases}
      d \, G (q/d; \, n/d, \, m/d) & \text{if} \quad d \mid m , \\
      0 & \text{if} \quad d \nmid m .
        \end{cases}
\]  
      
(ii) 
If $(q, 2n) = 1 $ then
\[
  G(q; n, m) = e_q \left( - \overline{(4n)} \, m^2 \right) \,
  \left( \frac{n}{q} \right) \, G(q; 1) .
\]

(iii) 
If $(q, 2)=1$ then
\[
    G^2(q; 1) = (-1)^{\frac{q-1}{2}} \, q .
\]
\end{lemma}
\hfill$\square$

In the next lemma we present A.Weil's estimate for the Kloosterman sum. 
For the proof we refer the reader to Chapter~11 of \cite{IwKo}.
\begin{lemma} \label{Kloosterman}
We have
\[
  |K(q; n, m)| \le \tau(q) \, q^{\frac{1}{2}} \, (q, n, m)^{\frac{1}{2}} .
\] 
\end{lemma}
\hfill$\square$

\bigskip

Next we establish the following 
\begin{lemma} \label{est-lambda}
If $ 8 \nmid q$ then
\begin{equation} \label{415}
 | \lambda(q; n, m) | \le 16 \, \tau^2(q) \, q^{\frac{1}{2}}  \, (q, n, m)^{\frac{1}{2}} .
\end{equation}
In particular we have
\begin{equation} \label{417}
\lambda(q) \ll q^{1 + \varepsilon} .
\end{equation}
\end{lemma}

\bigskip

{\bf Remark:} An estimate of type \eqref{415} holds for any positive integer $q$. We impose the condition $8 \nmid q$ 
because in this case the proof is slightly simpler and because in our analysis only such $q$ appear.

\bigskip

{\bf Proof:}
First we consider the case $2 \nmid q$.
It is clear that from \eqref{100}, \eqref{300} and \eqref{400} it follows
\begin{align}
  \lambda(q; n, m)
    & =
    \sum_{1 \le x, y \le q }
    e_q(n x + m y) \; q^{-1}
    \sum_{1 \le h \le q } 
    e_q(h (x^2 + y^2 + 1 ))
     \notag
     \\
     & = 
      q^{-1}
    \sum_{1 \le h \le q } 
      e_q(h) \; G(q; h, n) \; G(q; h, m) .
      \notag
      \\
        & = 
      q^{-1}
      \sum_{l \mid q}
    \sum_{\substack{ 1 \le h \le q  \\ (h, q)= \frac{q}{l} }} 
      e_q(h) \; G(q; h, n) \; G(q; h, m) .
      \notag
\end{align}
Now we apply Lemma~\ref{Gauss}, our assumption $2 \nmid q$ and the definition \eqref{200}. We get 
\begin{align}
  \lambda(q; n, m)
    & =
    q \sum_{\substack{l \mid q  \\  \frac{q}{l} \mid (m,n)   }} l^{-2}
    \sum_{\substack{1 \le r \le l \\ (r, l) = 1 }} e_l(r) \, G \left(l; r, nlq^{-1} \right) \, G \left(l; r, mlq^{-1}  \right)
    \notag \\
    & =
    q \sum_{\substack{l \mid q \\ \frac{q}{l} \mid (m,n) }} l^{-2} \, G^2(l, 1) \,
    \sum_{\substack{1 \le r \le l \\ (r, l) = 1 }} e_l \left( r - \overline{(4r)} (n^2 + m^2 ) l^2 q^{-2} \right) 
      \notag \\
    & =
    q \sum_{\substack{l \mid q \\ \frac{q}{l} \mid (m,n) }}  (-1)^{\frac{l-1}{2}} \; l^{-1} \;
    K \left( l; 1, \overline{4} (n^2 + m^2 ) l^2 q^{-2} \right) 
    \notag .
\end{align}
From the last formula and Lemma~\ref{Kloosterman} we find that if $2 \nmid q$ then
\begin{equation} \label{418}
  |\lambda(q; n, m)|
    \le
      q \sum_{\substack{l \mid q \\ \frac{q}{l} \mid (m,n) }} \tau(l) \, l^{-\frac{1}{2}} 
      \le q \, \tau(q) 
      \sum_{r \mid (q, n, m)} q^{-\frac{1}{2}} \, r^{\frac{1}{2}}
      \le \tau^2(q) \; q^{\frac{1}{2}} (q, n, m)^{\frac{1}{2}} .
\end{equation}

Next we note that the function $\lambda(q; n, m)$ possesses a multiplicative property.
More precisely, one can establish that if $(q_1, q_2)=1$ then
\begin{equation} \label{420}
  \lambda(q_1 q_2; n, m) = \lambda \left( q_1; n \overline{(q_2)}_{q_1}, m \overline{(q_2)}_{q_1} \right) \; 
  \lambda \left( q_2; n \overline{(q_1)}_{q_2}, m \overline{(q_1)}_{q_2} \right).
\end{equation}
(The proof is elementary and we leave it to the reader.) 
Now, in the general case, we represent $q = 2^{h} q_1$, where $2 \nmid q_1$ and $h \le 2$. We apply 
\eqref{418}, \eqref{420} and the trivial estimate $|\lambda(2^h; n, m)| \le 2^{2h}$ and obtain \eqref{415}.

\bigskip

Finally, it is clear that \eqref{417} follows directly from \eqref{415}.

\hfill
$\square$

\begin{lemma} \label{UV}
Suppose that $8 \nmid q$ and $D \ge 2$. Then for the sums
\[
   U =  \sum_{ 1 \le n \le D} \frac{|\lambda(q; n, 0)|}{n}  ,
     \qquad
   V =  \sum_{ 1 \le n, m \le D} \frac{|\lambda(q; n, m)|}{n \, m} .
\]
we have
\begin{equation} \label{500}
  U \ll q^{\frac{1}{2}+ \varepsilon} D^{\varepsilon} ,
  \qquad
  V \ll q^{\frac{1}{2}+ \varepsilon} D^{\varepsilon}.
\end{equation}
\end{lemma}

{\bf Proof:}
From Lemma~\ref{est-lambda} it follows that
\[
  U \ll q^{\frac{1}{2} + \varepsilon} \sum_{1 \le n \le D} \frac{(q, n)^{\frac{1}{2}}}{n}
    = q^{\frac{1}{2} + \varepsilon} \, \Sigma_0,
\]
say. Obviously
\begin{equation} \label{600}
  \Sigma_0 \ll \sum_{ r \mid q } r^{\frac{1}{2}} \sum_{\substack{ n \le D \\ n \equiv 0 \pmod{r}} } \frac{1}{n}
  \ll \sum_{ r \mid q} r^{- \frac{1}{2}} \, \log D \ll (q D)^{\varepsilon}
\end{equation}
and the first inequality in \eqref{500} follows.  To prove the second one we apply again Lemma~\ref{est-lambda} 
and use \eqref{600} to get
\[
  V \ll q^{\frac{1}{2} + \varepsilon} \sum_{1 \le n , m \le D} \frac{(q, n, m)^{\frac{1}{2}}}{n \, m}
  \ll q^{\frac{1}{2} + \varepsilon} \sum_{1 \le n , m \le D} \frac{(q, n)^{\frac{1}{2}} (q, m)^{\frac{1}{2}}}{n \, m}
  = q^{\frac{1}{2} + \varepsilon} \, \Sigma_0^2 \ll q^{\frac{1}{2} + \varepsilon} D^{\varepsilon} .
\]
\hfill
$\square$

\begin{lemma} \label{rho}
For any  $D \ge 2$ we have
\begin{equation} \label{1000}
  \rho(t) =  \sum_{1 \le |n| \le D} \frac{e(nt)}{2 \pi i n}
  + O \left( g(D, t)\right) , 
\end{equation}
where $g(D, t)$ is a positive, infinitely many times 
differentiable and periodic with period one 
function of $t$. It can be represented as a Fourier series
\[
  g(D, t) = \sum_{n \in \mathbb Z} c_D(n) \, e(nt) ,
\]
with coefficients $c_D(n)$ satisfying
\[
  c_D(n) \ll \frac{\log D}{D} \quad \text{for all} \quad n,
\]
and
\[
    \sum_{|n| > D^{1 + \varepsilon}} |c_D(n)| \ll D^{-A} .
\]
Here $A>0$ is arbitrarily large and the constant in the $\ll$--symbol depends on $A$ and $\varepsilon$.
\end{lemma}

{\bf Proof:}
It is well-known that $\rho(t)$ can be represented in a form similar to 
\eqref{1000}, but with a different form of the error term --- with the function 
$\min \left( 1, \left( D ||t|| \right)^{-1} \right) $ rather than $g(D, t)$. 
(For a proof we refer the reader to Chapter~2 of \cite{Hooley}.) 
The proof of the present formula for $\rho(t)$ is available in \cite{Tol}.

\hfill
$\square$

\section{Proof of the theorem}

\indent

We use the well-known identity $\mu^2(n) = \sum_{d^2 \mid n} \mu(d)$ to write
\[
  S(H) = \sum_{1 \le d \le \sqrt{2 H^2 + 1}} \mu(d) \, T \left(H, d^2 \right),
\]
where $T(H, q)$ is the number of pairs of positive integers $x, y \le H$
satisfying the congruence $ x^2 + y^2 + 1 \equiv 0 \pmod{q}$.
We choose a parameter $z$ satisfying
\begin{equation} \label{1050}
  \sqrt{H} \le  z \le H
\end{equation}
and denote by $S'$ the contribution to $S(H)$ coming from the terms with $z < d \le \sqrt{2H^2+1}$. 
It is clear that 
\[
  T \left( H, d^2 \right) = \sum_{1 \le l \le (2H^2 +1) \, d^{-2}} \;
   \sum_{\substack{ 1 \le x, y \le H \\ x^2 + y^2 = l d^2 -1}} 1 
  \ll H^{2 + \varepsilon} d^{-2} .
\]
Hence we obtain $ S' \ll H^{2 + \varepsilon} z^{-1}$ and therefore
\begin{equation} \label{1100}
  S(H) = \sum_{ 1 \le d \le z} \mu(d) \, T \left(H, d^2 \right) + O \left( H^{2 + \varepsilon} z^{-1}\right) .
\end{equation}

\bigskip

From this point onwards we assume that $q = d^2$,
where $d$ is squarefree and $d \le z$.
Denote by $ M(H, q, x ) $ the number of positive integers $h \le H$ satisfying
$h \equiv x \pmod{q}$.
Obviously
\begin{equation} \label{1150}
 M(H, q, x) = H \, q^{-1} + O \left( 1 \right) .
\end{equation}
It is clear that
\begin{equation} \label{1130}
  T \left( H, q \right) = \sum_{x, y : \, \eqref{400} } \,
    M(H, q, x) \, M(H, q, y) .
\end{equation}
If we use \eqref{1100} -- \eqref{1130} and choose $z = \sqrt{H}$ then
we may find an asymptotic formula for 
$S(H)$ with an error term $O \big( H^{\frac{3}{2} + \varepsilon} \big) $.
To establish our sharper result we represent the error term in \eqref{1150} in an explicit form.
First we write
\begin{equation} \label{1160}
   M(H, q, y ) =
  \left[ \frac{H-y}{q} \right] - \left[ \frac{-y}{q} \right] =
 \frac{H}{q} + \rho \left( \frac{H-y}{q} \right) - \rho \left( \frac{-y}{q} \right) .
\end{equation}
We substitute the last expression for $M(H, q, y)$ in \eqref{1130} 
and denote by $T'$ the contribution to $T(H, q)$ coming from the last term 
in the right-hand side of \eqref{1160}.
Next we represent $T'$ as
\begin{equation} \label{1200}
 T' = T^{\prime\prime} +  T^{\prime\prime\prime} ,
\end{equation}
where $ T^{\prime\prime} $ consists of the terms for which $x^2 + 1 \equiv 0 \pmod{q}$
and $ T^{\prime\prime\prime} $ is the contribution of the other terms.
We have
\begin{equation} \label{1250}
 T^{\prime\prime\prime} = 
 -
 \sum_{\substack{1 \le x \le q \\ x^2 + 1 \not\equiv 0 \pmod{q}}}
  M(H, q, x) 
 \sum_{\substack{1 \le y \le q \\ y^2 \equiv - x^2 - 1 \pmod{q}}}
   \rho \left( \frac{-y}{q} \right) = 0 
\end{equation}
because the last sum over $y$ vanishes. Indeed, it does not contain terms with
$y = \frac{q}{2}$ and $y = q$. Further, for any $y$
satisfying the congruence condition and such that $1 \le y < \frac{q}{2}$ 
the integer $q-y$ satisfies the same congruence and we have
$\rho \left( \frac{-y}{q} \right) + \rho \left( \frac{-(q-y)}{q} \right) = 0$.

\bigskip

Consider now $ T^{\prime\prime}$. 
We have
\begin{equation} \label{1270}
   T^{\prime\prime}  = 
 -
 \sum_{\substack{1 \le x \le q \\ x^2 + 1 \equiv 0 \pmod{q}}}
   M(H, q, x) 
 \sum_{\substack{1 \le y \le q \\ y^2 \equiv 0 \pmod{q}}}
   \rho \left( \frac{-y}{q} \right) 
   \ll H^{\varepsilon} \left( H q^{ - 1 } + 1 \right) .
\end{equation}
Indeed, the last sum over $y$ equals $O(1)$ because, 
according to the above arguments, it reduces to a sum with at most
two terms (corresponding to $y = \frac{q}{2}$ and $y = q$). Further, the number of solutions of the congruence
$x^2 + 1 \equiv 0 \pmod{q}$ equals $O \left( q^{\varepsilon} \right) $ and applying \eqref{1150}
we obtain \eqref{1270}.

\bigskip

From \eqref{1200} -- \eqref{1270} we get $T' \ll  H^{\varepsilon} \left( H q^{ - 1 } + 1 \right)$
and therefore
\[
  T \left( H, q \right) = \sum_{x, y : \, \eqref{400} }
   M(H, q, x) 
   \left( \frac{H}{q} + \rho \left( \frac{H-y}{q} \right) \right) + O \left(  H^{\varepsilon} \left( H q^{ - 1 } + 1 \right) \right) .
\]

We proceed with the quantity $M(H, q, x)$ in a similar manner and we conclude that
\[
  T \left( H, q \right) = \sum_{ x, y : \, \eqref{400} }
  \left( \frac{H}{q} + \rho \left( \frac{H-x}{q} \right) \right)
  \left( \frac{H}{q} + \rho \left( \frac{H-y}{q} \right) \right) + 
  O \left(  H^{\varepsilon} \left( H q^{ - 1 } + 1 \right) \right)  .
\]
Now we use \eqref{300} and \eqref{405} to find 
\begin{equation} \label{1290}
  T \left( H, q \right) = 
  \frac{H^2 \lambda(q)}{q^2} + 2\frac{H}{q} T_1 (H,q) + T_2 (H,q) +
   O \left(  H^{\varepsilon} \left( H q^{ - 1 } + 1 \right) \right)  ,
\end{equation}
where
\begin{align}
  T_1 (H,q)
    & =
     \sum_{ x, y : \, \eqref{400} }
     \rho \left( \frac{H-x}{q} \right) ,
     \label{1300} \\
   T_2 (H,q)
    & =
     \sum_{ x, y : \, \eqref{400} }  
     \rho \left( \frac{H-x}{q} \right)  \rho \left( \frac{H-y}{q} \right) .
     \label{1400} 
\end{align}

Consider $T_1(H, q)$. We apply \eqref{300} and Lemma~\ref{rho} with $D=H$ to get
\begin{equation} \label{2000}
  T_1(H, q) 
   =      T'_1(H, q) + O \left( T^*_1(H, q) \right) ,
\end{equation}
where
\begin{align}
  T'_1(H, q) 
  & =
  \sum_{ x, y : \, \eqref{400} } \;
   \sum_{1 \le |n| \le H}   \frac{e_q \left( n (H-x) \right)}{2 \pi i n }
   = 
   \sum_{1 \le |n| \le H}   \frac{e_q(nH) \, \lambda(q; -n, 0) }{2 \pi i n }
    ,
  \notag \\
  T^*_1(H, q)
  & =
    \sum_{ x, y : \, \eqref{400} } g\left( H, \frac{H-x}{q} \right) .
    \label{2050}
\end{align}

From Lemma~\ref{UV} it follows that 
\begin{equation} \label{2100}
  T'_1(H, q) \ll 
   H^{\varepsilon} q^{\frac{1}{2} } .
\end{equation}
To estimate $  T^*_1(H, q)$ we 
apply 
\eqref{300} and Lemmas \ref{est-lambda}, \ref{UV} and \ref{rho} and obtain
\begin{align}
  T^*_1(H, q) 
  & = 
   \sum_{ x, y : \, \eqref{400} } 
   \left( c_H(0) + 
   \sum_{1 \le |n| \le H^{1+ \varepsilon}} c_H(n) \, e_q(n(H-x))
   \right) + O (1)
   \notag \\
   & =
     c_H(0) \lambda (q) + 
     \sum_{1 \le |n| \le H^{1+ \varepsilon}} 
      c_H(n) \, e_q(nH) \, \lambda(q; -n, 0)
      + O (1)
      \notag \\
    & \ll 
      H^{\varepsilon - 1} q + 1 + 
      H^{\varepsilon - 1} 
      \sum_{1 \le |n| \le H^{1+ \varepsilon}} 
      |\lambda(q; -n, 0)|
    \notag  \\
    & \ll 
      H^{\varepsilon - 1} q + 1 + 
      H^{\varepsilon} 
      \sum_{1 \le n \le H^{1+ \varepsilon}} 
      \frac{|\lambda(q; n, 0)|}{n}
    \notag  \\
    & \ll 
      H^{\varepsilon - 1} q + H^{\varepsilon} q^{\frac{1}{2} }  .
      \label{2200}
\end{align}
From \eqref{2000}, \eqref{2100} and \eqref{2200} we get
\begin{equation} \label{2300}
  T_1(H, q) \ll H^{\varepsilon - 1} \, q +  H^{\varepsilon} \, q^{\frac{1}{2} }.
\end{equation}

Consider now $T_2(H, q)$. We apply \eqref{1400}, \eqref{2050}, \eqref{2200} and Lemmas~\ref{UV} and \ref{rho}, to get
\begin{align}
  T_2(H, q)
    & =
       \sum_{ x, y : \, \eqref{400} } \;
       \sum_{1 \le |n|, |m| \le H} \frac{e_q \left( (n+m) H \right)  \, e_q \left( - nx - m y \right) }
                                        { (2 \pi i )^2 n m }
                                        + O \left( H^{\varepsilon} T_1^* (H, q) \right) 
                                        \notag \\
    & =
     \sum_{1 \le |n|, |m| \le H}
       \frac{e_q \left( (n+m) H \right)  }
             { (2 \pi i )^2 n m } \,
               \lambda(q; -n, -m)
                 + O \left( H^{\varepsilon - 1} q + H^{\varepsilon} q^{\frac{1}{2} } \right) 
                                        \notag \\
    & \ll      
       \sum_{1 \le |n|, |m| \le H} \frac{|\lambda(q; n, m)|}{|nm|}
         +  H^{\varepsilon - 1} q + H^{\varepsilon} q^{\frac{1}{2} }
         \notag \\
    & \ll                                       
         H^{\varepsilon - 1} q + H^{\varepsilon} q^{\frac{1}{2} }  .
    \label{2400}
\end{align}

From \eqref{1290}, \eqref{2300} and \eqref{2400} it follows that
\[
   T(H, q) = H ^2 \frac{\lambda(q)}{q^2} + 
   O \left( H^{\varepsilon} (H q^{- \frac{1}{2} } + q^{\frac{1}{2}} + H^{-1} q)  \right) .
\]
We use the above formula, \eqref{1050} and \eqref{1100} and we get
\begin{align}
  S(H) 
    & = 
    H^2 \sum_{1 \le d \le z} \frac{\mu(d) \lambda \left( d^2 \right) }{ d^4 }
   + O \left( H^{\varepsilon} \left(  H + z^2 + H^{-1} z^3 + H^2 z^{-1}  \right)  \right)  
   \notag \\
    &  = 
    c \, H^2 + O \left( H^{\varepsilon} \left(H^2 z^{-1} + z^2 \right)  \right) ,
    \notag
\end{align}
where $ c = \sum_{d=1}^{\infty} \frac{\mu(d) \lambda \left( d^2 \right) }{d^4} $.
(Clearly this number coincides with the product at \eqref{412}.)
It remains to choose $z = H^{\frac{2}{3}}$ and the proof of the theorem is complete

\hfill
$\square$

\bigskip

{\bf Acknowledgments:}
The present research of the author is supported by Sofia University, Grant~172/2010.

\bigskip
\bigskip

\vbox{
\hbox{Faculty of Mathematics and Informatics}
\hbox{Sofia University ``St. Kl. Ohridsky''}
\hbox{5 J.Bourchier, 1164 Sofia, Bulgaria}
\hbox{ }
\hbox{Email: dtolev@fmi.uni-sofia.bg}}


\begin{thebibliography}{99}

\bibitem{Est1}
T.Estermann,
{\it Einige S\"atze \"uber quadratfreie Zahlen},
Math. Ann. 105, (1931), 653--662.

\bibitem{Est2}
T.Estermann,
{\it A new application of the Hardy--Littlewood--Kloosterman method},
Proc. London Math. Soc., 12 (1962), 425--444.

\bibitem{Filaseta}
M.Filaseta, 
{\it Powerfree values of binary forms}, 
J. Number Theory, 49 (2) 1994, 250-–268, 


\bibitem{Greaves}
G.Greaves, Power-free values of binary forms, Quart. J. Math. 
43 (1992), 45-–65.


\bibitem{Heath-Brown}
D.R.Heath-Brown,
{\it The square sieve and consecutive square-free numbers},
Math. Ann. 266 (1984), 251--259.

\bibitem{Hooley}
C.Hooley, 
{\it Applications of sieve methods to the theory of numbers},
Cambridge Univ. Press, 1976.


\bibitem{Hooley1}
C.Hooley, 
{\it On the power-free values of polynomials in two variables},
Roth 80th birthday volume, in press.

\bibitem{Hooley2}
C.Hooley, 
{\it On the power-free values of polynomials in two variables: II},
J. Number Theory 129 (2009) 1443–-1455


\bibitem{Hua}
L.-K.Hua, 
{\it Introduction to number theory},
Springer, 1982.

\bibitem{IwKo}
H.Iwaniec, E.Kowalski,
{\it Analytic number theory},
Colloquium Publications, vol. 53, Amer. Math. Soc., 2004.

\bibitem{Poon}
B.Poonen,
{\it Squarefree values of multivariable polynomials},
Duke Math. J., 118, 2 (2003), 353-373. 


\bibitem{Tol}
D.Tolev,
On the exponential sum with squarefree numbers,
Bull. London Math. Soc., 37, 6, (2005), 827-834. 


\end{thebibliography}
\end{document}